\documentclass[aps,pre,preprint,superscriptaddress]{revtex4-1}
\usepackage{graphicx}% Include figure files
\usepackage{epsfig}%add
\usepackage{subfigure}  %added
\usepackage{float}
\usepackage{braket}
\usepackage{pslatex,latexsym}
\usepackage{bm,float,lipsum}
\usepackage{amsmath}
\usepackage{amssymb}
\usepackage{amsfonts}
\usepackage{mathrsfs}
\usepackage{color}
\usepackage{longtable}
\usepackage{rotating}
\begin{document}

\title{Irrelevance of linear controllability to nonlinear dynamical networks}

\author{Junjie Jiang}
\affiliation{School of Electrical, Computer and Energy Engineering, Arizona State University, Tempe, Arizona 85287, USA}

\author{Ying-Cheng Lai}
\affiliation{School of Electrical, Computer and Energy Engineering, Arizona State University, Tempe, Arizona 85287, USA}
\affiliation{Department of Physics, Arizona State University, Tempe, Arizona 85287, USA}

\begin{abstract}

There has been tremendous development of linear controllability of complex networks. Real-world systems are fundamentally nonlinear. Is linear controllability relevant to nonlinear dynamical networks? We identify a common trait underlying both types of control: the nodal “importance”. For nonlinear and linear control, the importance is determined, respectively, by physical/biological considerations and the probability for a node to be in the minimum driver set. We study empirical mutualistic networks and a gene regulatory network, for which the nonlinear nodal importance can be quantified by the ability of individual nodes to restore the system from the aftermath of a tipping-point transition. We find that the nodal importance ranking for nonlinear and linear control exhibits opposite trends: for the former large-degree nodes are more important but for the latter, the importance scale is tilted towards the small-degree nodes, suggesting strongly irrelevance of linear controllability to these systems. The recent claim of successful application of linear controllability to {\em C. elegans} connectome is examined and discussed.

\end{abstract}

\maketitle

\section*{Introduction}

In the development of a field that involves dynamical systems, when knowledge
has accumulated to certain degree, the question of control would arise
naturally. For example, in nonlinear dynamics, the principle of controlling 
chaos was articulated in 1990~\cite{OGY:1990}, after approximately 
a decade of intense research focusing on the fundamental understanding of 
chaotic dynamical systems. Likewise, in complex networks, the issue of control 
began to be addressed~\cite{LH:2007,RJME:2009} also approximately 
after ten years of tremendous growth of research triggered by the pioneering 
work on small world and scale-free networks. A 
key development is the systematic adoption of the linear structural 
controllability theory to complex networks with directed 
interactions~\cite{LSB:2011}. Since then, there has been a 
great deal of effort in investigating the linear controllability of complex 
networks~\cite{WNLG:2012,NA:2012,YRLLL:2012,NV:2012,YZDWL:2013,MDB:2014,RR:2014,Wuchty:2014,YZWDL:2014,WBSS:2015,NA:2015,SCL:2015,LGS:2015,CWWL:2016,WCWL:2017,KSS:2017a,KSS:2017b}.

Control of linear dynamical systems is a traditional field in 
engineering~\cite{Kalman:1963,Lin:1974}. Because of the simplicity in the 
possible dynamical behaviors that a linear dynamical system can generate 
(in contrast to nonlinear dynamical systems where the behaviors are 
extremely rich and diverse), the general objective is to design proper 
control signals to drive the system from an arbitrarily initial state to 
an arbitrarily final state in finite time. When applying the linear 
controllability theory to complex networks, a primary goal has been to 
determine the minimum number of controllers. This problem was 
addressed~\cite{LSB:2011} for complex directed networks through the 
development of a minimum input theory based on the concept of maximum 
matching~\cite{HK:1973,ZOY:2003,ZM:2006}. To generalize the linear 
controllability theory to networks of arbitrary structures (e.g., weighted 
or unweighted, directed or undirected), an exact controllability theory was 
developed~\cite{YZDWL:2013} based on the Popov-Belevitch-Hautus (PBH) rank 
condition~\cite{PBH:1969}. The exact controllability 
theory provides a computationally extremely efficient method to determine not 
only the minimum number of controllers but also the set of the nodes to which 
the control signals should be applied - the set of driver nodes, for complex 
networks of arbitrary topology and link structures~\cite{YZDWL:2013}.   

The development of the linear controllability theories has played the role
of stimulating research on controlling complex networks~\cite{LB:2016}. 
However, its limitations must not be forgotten. The fundamental assumption 
used in any linear controllability theory is that the nodal dynamics are 
described by a set of coupled linear, first-order differential equations. 
While such a setting may be relevant to engineering control systems, 
real-world systems are governed by nonlinear dynamics such as biologically 
inspired networks~\cite{ABM:2018}. In classical control engineering, it is 
well recognized that controllability for nonlinear systems requires a 
different set of tools to be developed compared to what is known for the 
controllability of linear systems~\cite{NV:book}.
A serious concern is the tendency to 
overstate the use or the predictive power of the linear controllability 
theories when they are applied to real-world physical or biological systems. 
For example, it was claimed recently~\cite{YVTCWSB:2017} that linear network 
control principles can predict the neuron function in the {\em Caenorhabditis 
elegans} connectome, a highly nonlinear dynamical neuronal network. 
The goal of the present work is to legitimize this concern in a quantitative
manner by presenting concrete and statistical evidence that linear network 
controllability may not be relevant to physically or biologically meaningful 
control of nonlinear networks. 

The physical world is nonlinear. Network dynamics in biological or ecological 
systems are governed by nonlinear rules with no exceptions. Control of real 
world complex networks based on the rules of nonlinear dynamics has remained 
to be an extremely difficult problem. Existing strategies include local 
pinning~\cite{WC:2002,LWC:2004,SBGC:2007,YCL:2009}, feedback vertex set 
control~\cite{FMKS:2013,MFKS:2013,ZYA:2017}, controlled switch among 
coexisting attractors~\cite{WSHWWGL:2016}, or local control~\cite{KSS:2017b}. 
These methods belong to the category of open-loop control, i.e., 
one applies pre-defined control signals or parameter perturbations to a 
feedback vertex set chosen according to some physical criteria. 
For certain nonlinear dynamical networks, especially those in ecology,
closed-loop control can be articulated and has been demonstrated to be 
effective~\cite{SLLGL:2017}. Recently, how to exploit 
biologically inspired agent based control method to choose different 
alternative states in engineered multiagent network systems has been 
studied~\cite{GFSL:2018}. 

In order to answer the question ``is linear controllability relevant to 
nonlinear dynamical networks?'', two challenges must be met. Firstly, 
because of lack of general controllability framework for nonlinear networks 
it is necessary to focus on {\em specific} contexts where nonlinear network 
control can be done in a physically or biologically meaningful way. We choose 
two such contexts: mutualistic networks in ecology~\cite{BJMO:2003,GJT:2011,
NJB:2013,LNSB:2014,RSB:2014,DB:2014,GPJBT:2017,JHSLGHL:2018} and a gene 
regulatory network from systems biology~\cite{Alon:2006,BBILA:2006,GBB:2016}.
Secondly and more importantly, linear and nonlinear dynamical networks are
fundamentally and characteristically different in many aspects, so are the 
respective control methods. How do we compare their control performances? 
(How can an apple be compared with a banana?) Our idea is that, even in the 
analog of apple-banana comparison, if one finds a common trait, e.g., the 
amount of sugar contained per gram of the substance, then a comparison between 
an apple and a banana in terms of the specific common trait is meaningful. 
We are thus led to seek a feature or a characteristic that is common in both 
nonlinear and linear network control. Specifically, we identify the statistical 
importance of individual nodes in control as such a common trait. 

Our approach and main result can be described, as follows. Given a nonlinear 
dynamical network with its structure determined from empirical data, we focus 
on the concrete problem of harnessing a tipping point at which the system 
transitions from a normal state to a catastrophic state (e.g., massive 
extinction) or from a catastrophic state to a normal state abruptly as a 
system parameter changes through a critical point~\cite{Schefferetal:2009,
Scheffer:2010,WH:2010,DJ:2010,DVKG:2012,BH:2013,TC:2014,LNSB:2014,JHSLGHL:2018}.
%In the former case where the tipping point transition is undesired, to 
%develop a biologically viable control strategy to remove the tipping point 
%so as to delay the occurrence of global extinction is of broad interest. 
%In the latter case where the transition can land the system in a normal 
%state, to induce a tipping point to achieve global restoration is desired. 
We exploit the ability of the individual nodes, via control, to make the 
system recover from the aftermath of a tipping point transition that puts 
the system in an extinction state. This enables a quantitative ranking of 
the importance of the individual nodes to be determined. The ranking is 
generally found to be linearly correlated with the nodal degree of the network, 
in agreement with intuition. The individual nodes, in 
terms of their ability to make the system recover, are drastically distinct.
We then perform linear control on the same network by assuming 
artificial linear nodal dynamics. Using the exact controllability 
theory~\cite{YZDWL:2013}, we calculate the minimal control set. 
A key feature of linear network control, which was usually not emphasized in
most existing literature on linear controllability~\cite{WNLG:2012,NA:2012,YRLLL:2012,NV:2012,YZDWL:2013,MDB:2014,RR:2014,Wuchty:2014,YZWDL:2014,WBSS:2015,NA:2015,SCL:2015,LGS:2015,CWWL:2016,WCWL:2017,KSS:2017a,KSS:2017b} but was
mentioned in a recent paper~\cite{CARRSA:2017}, is that the minimal 
control set of nodes is not unique. 
For a reasonably large network (e.g., of size of a few hundred), there can be
vastly many such sets that are equivalent to each other in terms of control
realization. Thus, in principle, there is a finite probability for a node in 
the network to be chosen as a control driver and the corresponding probability 
can be calculated from the ensemble of the minimal control sets. This 
probability can be defined as a kind of importance of the node in control 
relative to other nodes so that a nodal importance ranking can be determined. 
Because of the generality and universality of the linear control framework, 
the method to determine the nodal importance is applicable to any complex 
network. For a large number of real pollinator-plant mutualistic networks 
reconstructed from empirical data from different geographical regions of the 
world (Table II) and a representative gene regulatory network, 
we find that the linear 
importance ranking favors the small degree nodes, in stark contrast to the 
case of nonlinear control where large degree nodes are typically more valuable.
The characteristic difference in the importance ranking of the nodes in terms
of their role in control, linear or nonlinear, suggests that linear   
controllability may not be relevant to physically or biologically justified 
nonlinear control for the mutualistic and gene regulatory networks.

\section*{Results}

\subsection*{A concrete example of complex pollinator-plant mutualistic network illustrating irrelevance of linear controllability}

The assumptions of this study are as follows. For linear dynamical networks, 
a general controllability framework exists, which can be used to determine
the nodal importance ranking and is applicable to all networks. For 
nonlinear networks, because of the rich diversity in their dynamics, at the 
present a general control framework does not exist. The control strategy thus 
depends on the specific physical or biological context of the network.

\begin{figure}[ht!]
\centering
\includegraphics[width=\linewidth]{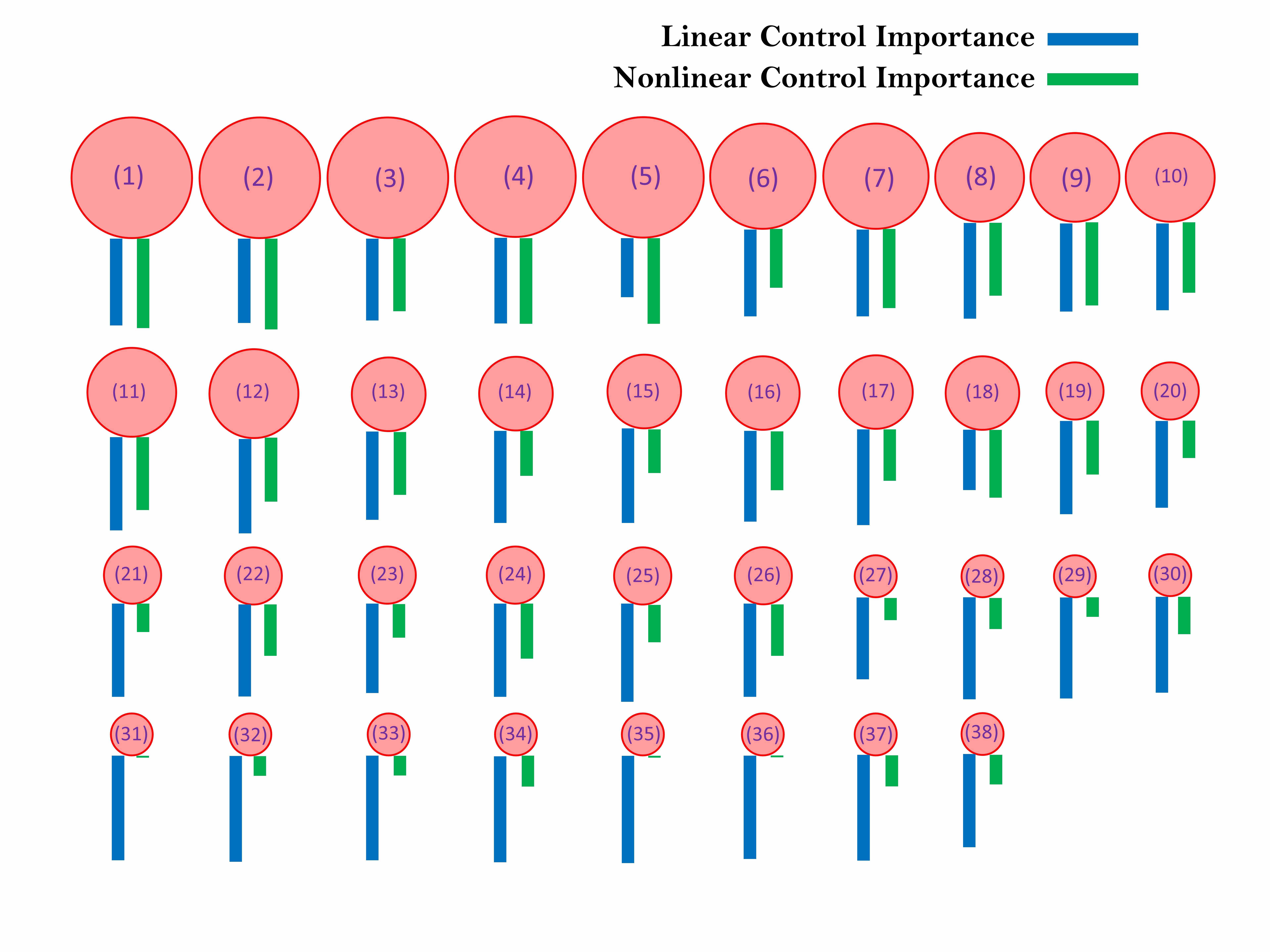}
\caption{ \textbf{Distinct characteristics in nonlinear 
and linear control of a representative complex mutualistic network.}
The system is network $A$ reconstructed based on empirical data from Tenerife, 
Canary Islands~\cite{DHO:2003}. The numbers of pollinators, plants, and 
mutualistic links are $N_A=38$, $N_P=11$, and $L=106$, respectively. For each 
node, the species name is given in Table I. The
length of the green bar below each species is indicative of the relative 
importance of the node in tipping point control of the actual nonlinear 
dynamical network, which is calculated based on Eq.~(\ref{eq:def_R_N}).
The blue bars illustrate the relative importance of the nodes when the system 
is artificially treated as a linear, time-invariant network, which are 
calculated according to Eq.~(\ref{eq:def_R_L}). There is great variation 
in the lengths of the green bars for different species, demonstrating a 
highly non-uniform nonlinear control importance ranking. In contrast, there is 
little variation in the length of the blue bars among the different species,
indicating an approximately uniform linear control importance ranking.
Linear controllability may thus not be useful for controlling the actual 
nonlinear dynamical network.}   
\label{fig:SFMN}
\end{figure}

To demonstrate the characteristic statistical difference between nonlinear
and linear control, we take a representative pollinator-plant mutualistic
network (network $A$), and calculate the node based, nonlinear and linear 
control importance according to Eqs.~(\ref{eq:def_R_N}) and 
(\ref{eq:def_R_L}), respectively, as described in {\bf Methods}. 
Figure~\ref{fig:SFMN} shows the 38 pollinator and plant 
species, together with the relative nonlinear and 
linear control importance as represented by the lengths of the green and blue 
bars beneath the images, respectively. There is a wide spread in the nonlinear 
control importance, but the linear control importance appears approximately 
uniform across the species. There are cases where a node is not important at 
all for nonlinear control (e.g., the first, fifth, and sixth species in the 
bottom row), but the node is important for linear control. The statistical 
characteristics of the nodal importance in nonlinear and linear control are 
thus drastically distinct. An examination of other empirical mutualistic 
systems reveals that, for some networks,
the behaviors are similar to those in Fig.~\ref{fig:SFMN}, while in others, 
the nodal importance shows opposite trends in nonlinear and linear control.
For example, there are cases where the nonlinear control importance tends 
to increase with the nodal degree, but the linear nodal importance shows 
the opposite trend. These results suggest that linear controllability may
not be useful for controlling the actual nonlinear dynamical network.

\subsection*{Nonlinear and linear control importance}

\begin{figure}[ht!]
\centering
\includegraphics[width=\linewidth]{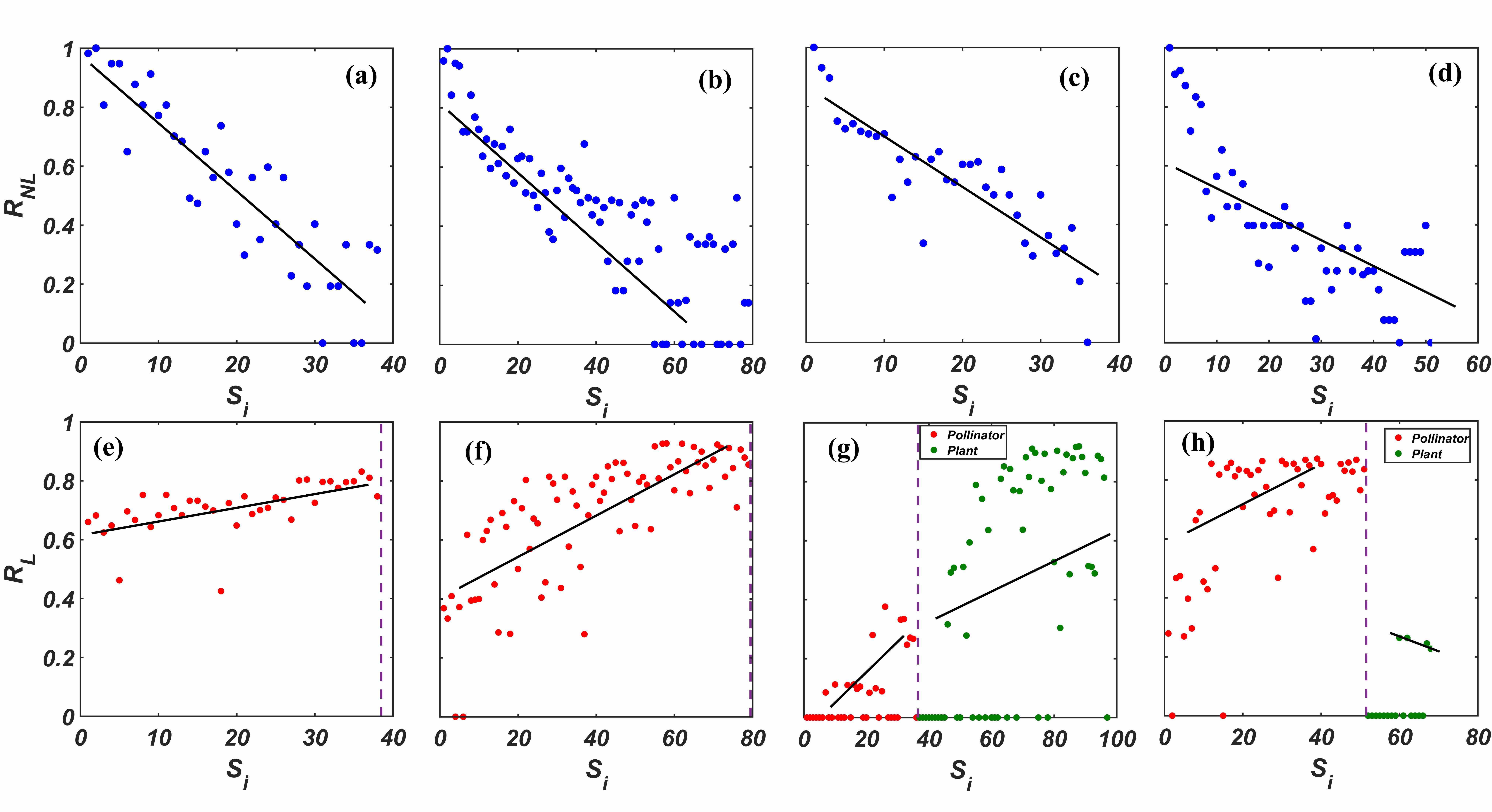}
\caption{ \textbf{Contrasting behaviors of nodal importance ranking in 
nonlinear and linear control}.
The four empirical networks are labeled as $A$, $B$, $C$, and $D$ with details 
given in {\bf Methods}. (a-d) Nonlinear and (e-h) linear control importance 
ranking for networks $A-D$, respectively. For tipping point control of the 
nonlinear network in (a-d), only the pollinator species are subject to 
external intervention through the managed maintenance of the abundance of 
a single species. 
The nodal index on the abscissa of each panel is arranged according to the 
degree ranking of the node: from high to low degree values (left to right). 
For the set of nodes with the same degree, their ranking is randomized.	
The nonlinear control importance is calculated from Eq.~(\ref{eq:def_R_N}) for 
the parameter setting $h=0.2$, $t=0.5$, $\beta_{ii}^{(A)}=\beta_{ii}^{(P)}=1$, 
$\beta_{ij}^{(A)}=\beta_{ij}^{(P)}=0,\alpha_i^{(A)} =\alpha_i^{(P)}=-0.3$,
and $\mu_A=\mu_P=0.0001$. The coupled nonlinear differential equations are
solved using the standard Runge-Kutta method with the time step $0.01$. 
The distinct feature associated with nonlinear control is that, in spite of
the fluctuations, larger degree nodes tend to be more important (i.e., 
more effective in recovering the species abundances after a tipping point). 
The linear control importance ranking in (e-h) can be calculated for all 
species based on definition (\ref{eq:def_R_L}), because the corresponding 
artificial linear dynamical network does not distinguish between pollinator 
and plant species. In each panel, the pollinators (red dots) and plants 
(green dots) are placed on the left and right side, respectively, and are 
arranged in descending values of their degree, with a vertical dashed line 
separating the two types of species. The striking result is that, for the 
pollinators, their ranking of linear control importance exhibits a trend 
opposite to that of nonlinear control importance.
A similar behavior occurs for ranking based on betweenness centrality
and actual degrees.}  
\label{fig:4Real_MN}
\end{figure}

\begin{figure}[ht!]
\centering
\includegraphics[width=0.8\linewidth]{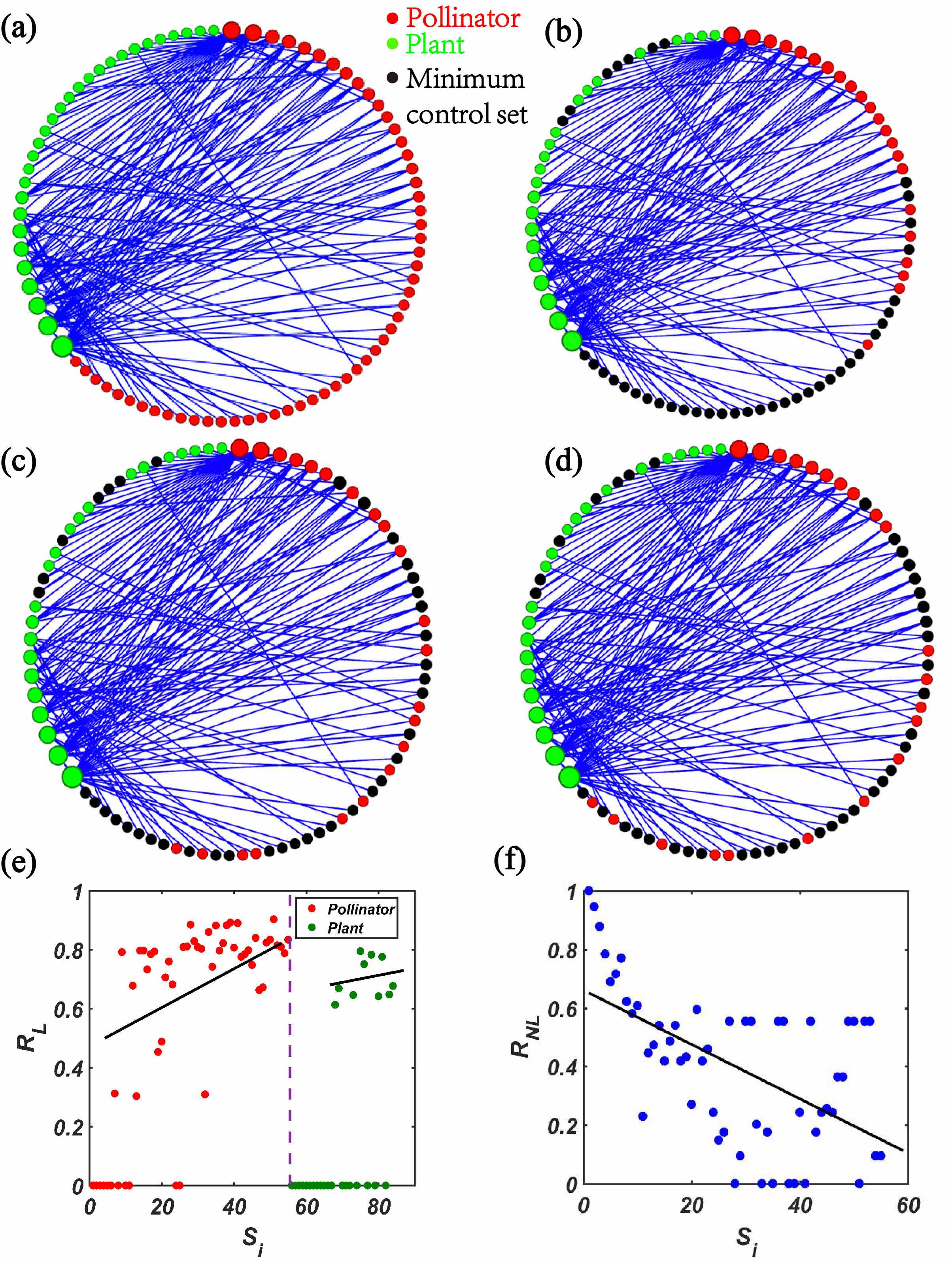}
\caption{ \textbf{Examples of distinct minimum controller sets associated with 
linear control}. For mutualistic network $E$ as described in {\bf Methods}, 
(a) network structure, where the size of a circle (red and green for a 
pollinator and a plant, respectively) is proportional to the degree of this 
node, (b-d) three examples of minimum controller sets (black dots), (e) linear 
control importance ranking, and (f) nonlinear control importance ranking. 
Other parameters are the same as those in Fig.~\ref{fig:4Real_MN}.}
\label{fig:LC_NC_Net_12}
\end{figure}

We present quantitative results of nonlinear and linear control importance for 
four empirical mutualistic networks described in {\bf Methods}, as shown in
Fig.~\ref{fig:4Real_MN}. For a given empirical network, to calculate the 
nonlinear control importance based on definition (\ref{eq:def_R_N}), we
begin from a zero value of the average mutualistic interaction strength
$\gamma_0$ where the system is in an extinction state without control, apply 
the control by setting the abundance of a pollinator species at 
$A_S=1.5$, and systematically increase the value of $\gamma_0$ towards a 
relatively large value (e.g., 3.0). During this process, the recovery point
$\gamma_c^i$ can be obtained. When the values of the recovery point for all
pollinator species have been calculated, Eq.~(\ref{eq:def_R_N}) gives the 
control importance for each species, as shown in Figs.~\ref{fig:4Real_MN}(a-d)
for networks $A-D$, respectively, where the index of the pollinator species
on the abscissa is arranged according to the nodal degree. Apart from 
statistical fluctuations, there is a high level of positive correlation 
between the nonlinear control importance and degree, i.e., larger degree
nodes tend to be more important. In particular, managed 
control of larger degree nodes is more effective for species recovery.    
To obtain the linear control importance according to Eq.~(\ref{eq:def_R_L}),
we use 1000 random minimum controller sets as determined by the linear exact
controllability to calculate the probability for each species to be chosen as
a driver node. Note that, because of the artificial imposition of linear time 
invariant dynamics on each node, there is a probability for any species to be 
a driver node, regardless of whether it is a pollinator or a plant species. 
The results are presented in Figs.~\ref{fig:4Real_MN}(e-h) for networks
$A-D$, respectively, where the linear control importance of the pollinators
(red dots) and that of the plants (green dots) - separated by the vertical
dashed line, are shown. The common feature among the four empirical networks
is that the linear control importance ranking has an opposite trend to the
nonlinear control importance ranking. That is, smaller degree nodes tend to
be more important for linear control. The correlation between linear control 
importance and degree is thus negative, which is in stark contrast to the 
behavior of nonlinear control importance. Overall, 
Figs.~\ref{fig:4Real_MN}(a-h) reveal that, for nonlinear control of tipping 
points, managing large degree nodes can be significantly more effective than 
harnessing small degree nodes, but for linear control of the same network, 
the large degree nodes play little role in control as they rarely appear 
in any minimum controller set.

The linear control importance measure, as defined in Eq.~(\ref{eq:def_R_L}), 
is rooted in the fact that, in the linear controllability theory, typically
there are many equivalent minimum controller sets~\cite{CARRSA:2017}.
It is useful to 
visualize such sets. Figure~\ref{fig:LC_NC_Net_12}(a) exhibits a graphical
representation of an empirical mutualistic network - network $E$ described
in {\bf Methods}, where the pollinators (red dots) and plants (green dots)
are arranged along a circle, and the size of a dot is proportional to the
degree of the corresponding node. By definition, mutualistic interactions
mean that there are no direct links between any pair of dots with the same 
color - any link in the network must be between a red and a green dot. For
this network, there are altogether approximately $10^{12}$ minimum controller 
sets of exactly the same size - three examples are shown in 
Figs.~\ref{fig:LC_NC_Net_12}(b-d), 
respectively, where the driver nodes are represented by black dots. A 
feature is that the minimum controller sets tend to avoid nodes of very 
large degrees in the network, which is consistent with the results in 
Fig.~\ref{fig:4Real_MN}. 
The corresponding linear and nonlinear control importance rankings are shown 
in Figs.~\ref{fig:LC_NC_Net_12}(e) and \ref{fig:LC_NC_Net_12}(f), respectively.
A comparison of these results indicates that the ranking behaviors are 
characteristically distinct, suggesting the difference 
between linear controllability and nonlinear control - the same 
message conveyed by Fig.~\ref{fig:4Real_MN}.

\subsection*{Gene regulatory networks}

\begin{figure}[ht!]
\centering
\includegraphics[width=\linewidth]{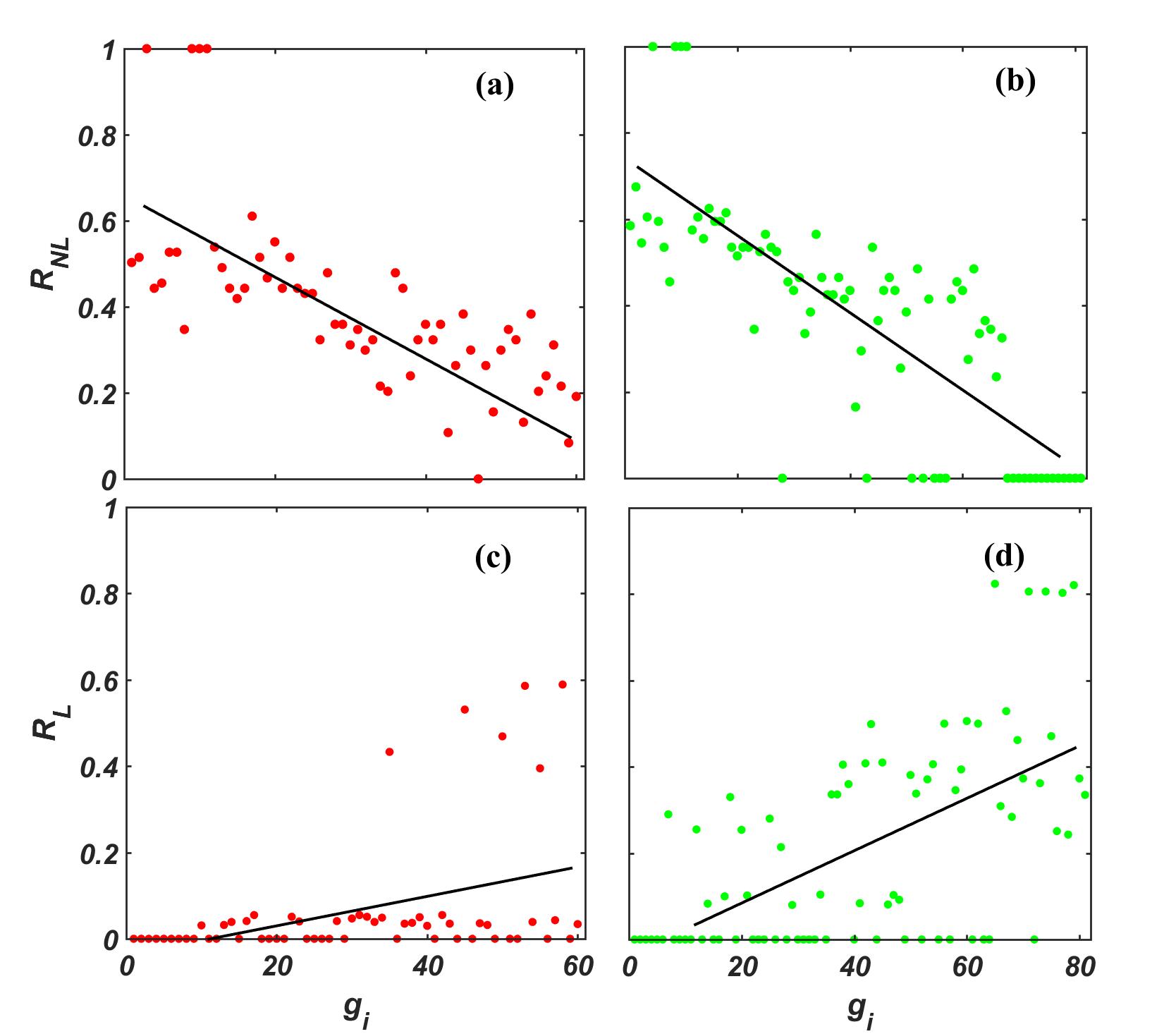}
\caption{ \textbf{Nodal importance rankings associated with nonlinear and 
linear control of a gene regulatory network.} (a,c) Nonlinear and linear 
control importance rankings for the subnetwork of the giant component of size 
$60$, respectively. In (a), there are four nodes with $R_{NL} = 1$, because 
controlling any of these genes will make the system recover immediately 
from the tipping point collapse when the direction of the change in the 
bifurcation parameter is reversed. For linear control in (c), the size of 
any minimum controller set is $N_D = 4$. (b,d) Nonlinear and linear nodal 
importance rankings for the subnetwork of $81$ nodes with input connection. 
In (b), there are several genes with $R_{NL} = 0$, as each gene in this 
group lacks the ability to restore the entire system even when its activity 
level is maintained at a high level through external control. For linear 
control in (d), any minimum controller set has $N_D = 17$ nodes. In all 
panels, the nodal index along the abscissa is arranged in the descending 
order of the outgoing degrees of the genes.}   
\label{fig:Gene_net}
\end{figure}

The opposite behaviors in the nodal importance ranking for linear 
controllability and nonlinear control also arise in gene regulatory networks. 
For such networks, tipping point dynamics similar to those in mutualistic 
networks can occur when a biological parameter is reduced, rendering 
feasible a similar control strategy (see {\bf Methods}). 
Figure~\ref{fig:Gene_net} shows, for the network of 
\textit{S. cerevisiae} described in {\bf Methods}, the nonlinear and linear 
control importance rankings for two subnetworks: the giant component (a,c) and 
the subnetwork of all nodes with input connections (b,d). Because of the dense 
connectivity in the giant component subnetwork, for linear control the size 
of the minimum controller set is $N_D = 4$ (c). For the subnetwork in (b,d), 
we have $N_D = 17$. Note that, for nonlinear control of the subnetwork (b), 
there are several genes that have zero nonlinear control importance, i.e., 
external management of the activation level of any of these genes is unable 
to restore the network function destroyed by a tipping point transition. The 
striking finding is that, for linear control, these genes are exceptionally 
important because the probability for any of these genes to belong to a 
minimum controller set is disproportionally high (e.g., $> 80\%$). 
If one follows the prediction of the linear controllability theory to      
identify those nodes as important and attempts to use them as the 
relevant nodes for actual control of the nonlinear network, one would be 
disappointed as harnessing any of these genes will have no effect on the 
tipping-point dynamics of the network. The occurrence of such genes with zero 
nonlinear control importance is the result of the interplay between the 
Holling-type of nonlinear dynamics and the complex network structure.

\subsection*{Pearson correlation and cosine distance}

\begin{figure}[ht!]
\centering
\includegraphics[width=0.9\linewidth]{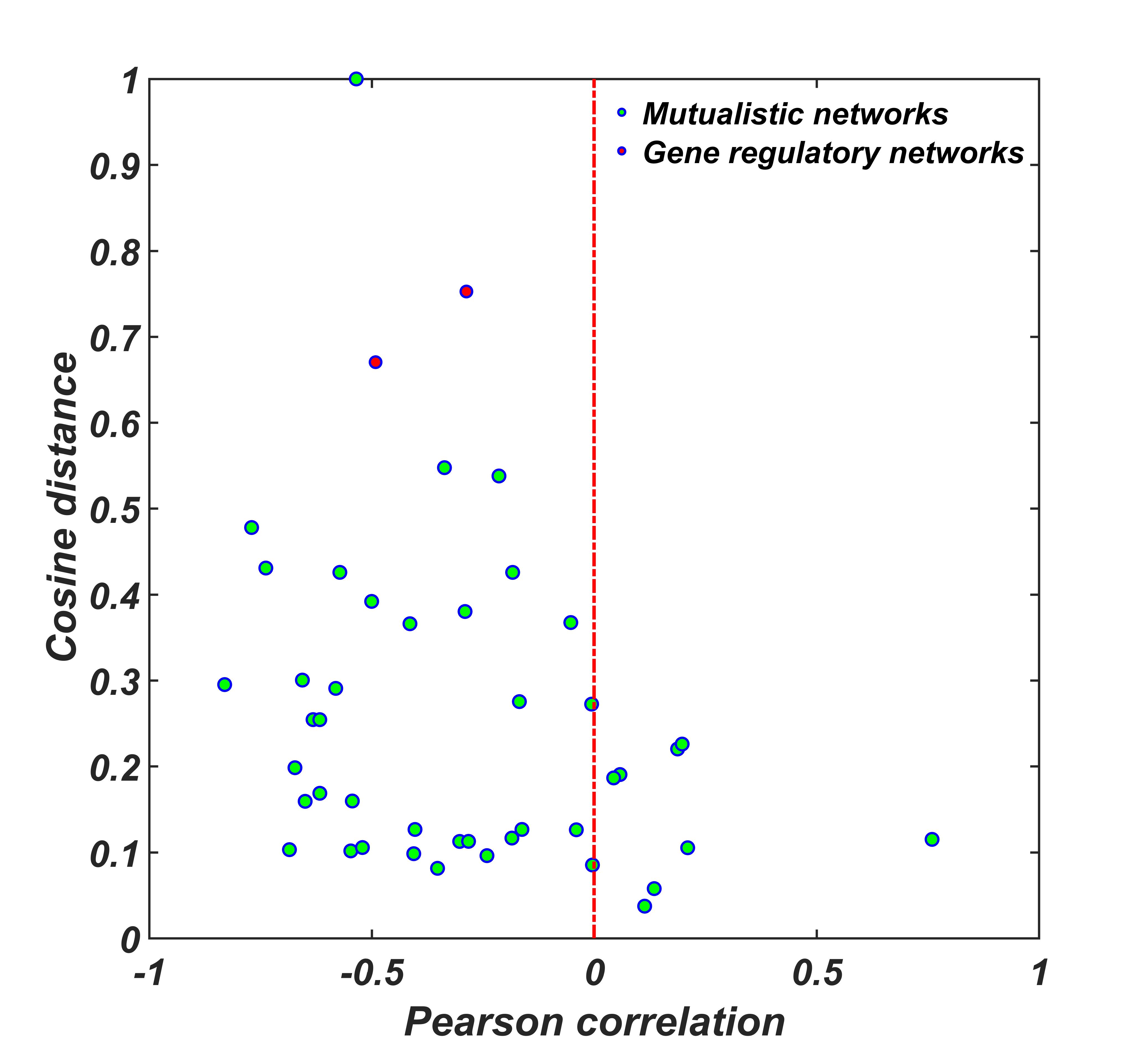}
\caption{ \textbf{Pearson correlation and cosine distance between linear 
and nonlinear control importance}. The abscissa and ordinate correspond to 
the values of Pearson correlation and cosine distance, respectively, 
between linear and nonlinear control importance. Each green circle corresponds
to a real mutualistic network (there are 43 of them) and the two red dots 
are for the two gene regulatory subnetworks in Fig.~\ref{fig:Gene_net}.
If there were a kind of relevance between nonlinear and linear control of
the same network, the dots would concentrate in the lower right region of
the plane with positive Pearson correlation and a small cosine distance. For
most of the empirical networks tested, the dots are in the region of negative
correlation with cosine distances below 0.4. For the two gene regulatory
subnetworks, not only are the values of the Pearson correlation negative, 
the cosine distances are also large.}  
\label{fig:CR_Cos_dis}
\end{figure}

For the five mutualistic networks ($A-E$) and two gene regulatory subnetworks
tested so far, the correlation between nonlinear and linear control importance 
is negative, as shown in Figs.~\ref{fig:4Real_MN}-\ref{fig:Gene_net}.
To test if this holds for a broad range of empirical networks, we calculate
the Pearson correlation and the cosine distance between linear and 
nonlinear control importance for a large number of real networks, as shown
in Fig.~\ref{fig:CR_Cos_dis}. In most cases, the correlation is negative and
the cosine distance is large. There are a few mutualistic networks with 
positive but small correlation. Out of the 43 mutualistic networks, only 
one has a large correlation value and a
small cosine distance (one corresponding to the rightmost green circle).
A peculiar feature of this network is that it has only six 
pollinator species and any minimum controller set in linear control contains 
four such species, rendering atypical this case.

Our detailed comparison between the control importance ranking in a type of 
biologically meaningful nonlinear control and in linear control for a large
number of real pollinator-plant mutualistic networks and a gene regulatory 
network provides evidence that linear controllability may generate results that
are drastically inconsistent with nonlinear dynamical behaviors and control of 
the system. In no way should this be a surprise, as the assumption of linear, 
time-invariant dynamics cannot be expected to hold for nonlinear dynamical 
networks in the real world. However, there is a recent tendency to apply the 
linear controllability framework to real-world nonlinear systems 
such as the \textit{C. elegans} connectome~\cite{YVTCWSB:2017} 
and brain networks\cite{GPCTAKMVMG:2015,MPGCGVB:2016,Tangetal:2017,TB:2018}. 
Although the linear control framework may provide insights into nonlinear 
dynamical networks under some specific circumstances, controlling highly 
nonlinear dynamical networks is still an open problem at the present. 
Nonetheless, a thorough analysis of the linear controllability would give 
clues to its inappropriateness and likely failure in real world systems 
(see Appendix A and Figs.~\ref{fig:LC_elegans} and \ref{fig:C_elegans_SPC}).

\section*{Discussion}

It is apparent that the assumption of linear, time invariant nodal dynamics 
is not compatible with natural systems in the real world that are governed by 
nonlinear dynamical processes. Why then study the linear controllability 
of complex networks? There were two reasons for this. Firstly, when the 
development of the field of complex network had reached the point at which 
the problem of control emerged as a forefront problem (around 2011), to 
adopt linear controllability, a well established framework in traditional 
control engineering, to complex networks seemed to be a natural starting 
point. The well developed mathematical foundation of linear control 
made it possible to address the effect of complex network structure on the
controllability in a rigorous manner~\cite{LSB:2011,YZDWL:2013}, physical
or biological irrelevance notwithstanding. Secondly,
to study the linear controllability of complex networks is justified from
the point of view of engineering, as linear dynamical systems are relevant 
to subfields in engineering such as control and signal processing. That 
being said, the applicability of the linear controllability to real physical,
chemical and biological systems is fundamentally limited because of the 
ubiquity of nonlinear dynamics in natural systems - a well accepted 
fact, thanks to more than four decades of extensive and intensive
study of nonlinear dynamics and chaos theory. It is imperative and a 
common sense understanding that the linear controllability of complex networks 
not be overemphasized and its importance and significance not be overstated.  

Quite contrary to the common sense understanding, there are recent 
claims that linear network controllability is applicable to real biological 
systems~\cite{YVTCWSB:2017,GPCTAKMVMG:2015,MPGCGVB:2016,Tangetal:2017,TB:2018}
for gaining new understanding. Curiosity demands a thorough reexamination 
of these claims. More importantly, such claims, if they are indeed unjustified 
but remain uncorrected, can potentially generate undesirable and negative 
impacts on the further development of the field of complex network control. 
These considerations motivated our present work.

The main question we have set out to answer is whether linear controllability 
is actually relevant to controlling nonlinear dynamical networks. To be able to
address this question, it is necessary to have nonlinear networked systems
for which a certain type of physically or biologically meaningful control 
can be carried out. We have identified two classes of such systems: complex 
pollinator-plant mutualistic networks in ecology and gene regulatory networks
in systems biology. We focus on the physically significant issue of 
controlling tipping points, which enables the nodal importance in the control
to be ranked. This is essentially a ranking associated with nonlinear control.
Ignoring the nonlinear dynamics and simply using the network structure to 
treat it as a linear, time-invariant system enable us to calculate the minimum
controller set in the linear controllability framework. Taking advantage
of the exact controllability theory~\cite{YZDWL:2013}, we identify a large 
number of equivalent configurations of the minimum controller set and find 
that, typically, there is a probability for almost every node to be in such a 
set. This probability serves as the base for ranking the nodal importance in 
linear controllability. The two types of control importance rankings, one 
nonlinear and another linear, can then be meaningfully compared. The main
finding of this paper is that the nonlinear and linear rankings are 
characteristically different for a large number of real world mutualistic
networks and the gene regulatory network of \textit{S. cerevisiae}. In 
particular, the nonlinear control importance ranking typically exhibits a 
behavior that in general favors high degree nodes. 
However, linear ranking typically exhibits the opposite trend that favors 
small degree nodes. These results are evidence that linear controllability 
theory generates information that is not useful for nonlinear control of
tipping point dynamics in complex biological networks.
A quite striking finding is that, for the gene regulatory network of 
\textit{S. cerevisiae}, there are four genes with essentially zero nonlinear 
control importance in the sense that managed control of any of these genes 
is unable to recover the system from the aftermath of a tipping point 
transition. However, in linear control, these four genes are far more 
important than other nodes in the network. Thus, for the particular gene 
regulatory network studied here, linear controllability absolutely has nothing 
to do with the actual control of the nonlinear dynamical network.   

In a recent work~\cite{YVTCWSB:2017}, it was claimed that linear structural 
controllability predicts neuron function in the {\em C. elegans} connectome. 
This real neuronal network has about 300 neurons, which contains four different
types of neurons including the sensory neurons, inter-neurons, and motor 
neurons. A sensory neuron can generate an action potential propagating to other
neurons, while an inter-neuron can receive action potentials from sensory 
neurons or other inter-neurons. The processes of generating and propagating 
action potentials are highly nonlinear. The claim of Ref.~\cite{YVTCWSB:2017} 
is thus questionable. We find that the {\em C. elegans} 
connectome, when artificially treated as a linear network, is uncontrollable
if the control signals are to be applied to sensory neurons only. A 
calculation of the linear control importance reveals an approximately uniform 
ranking across all neurons. The surprising feature is that, on average, a 
muscle cell is almost twice as important as a motor neuron in terms of 
linear controllability, but biologically any control signal must flow from 
neurons to muscle cells, not in the opposite direction. Linear controllability 
thus yields a result that is apparently biologically meaningless. In fact,
the ability to predict neuron function is based on signal 
propagation from some sensory to some motor neurons, which can be 
accomplished through random stimulation of some sensory neurons. Because of 
the existence of great many equivalent minimal control driver sets, which
sensory neuron should be chosen to deliver a control signal is completely 
random. From the point of view of signal paths, there exist vastly large
numbers of direct paths from the sensory to the motor neurons. Because of the
approximately uniform ranking in nodal importance as a result of the existence
of many equivalent minimum controller sets, linear controllability theory,
when being used fairly in the sense of taking into considerations of the many
controller set realizations, cannot possibly yield any path that is more 
special than others to uncover hidden biological functions. 
That is, it is not necessary to use linear controllability to predict 
any neuron function, contradicting the claim in Ref.~\cite{YVTCWSB:2017}.
If control were to play a role in predicting some functions, it must be some 
kind of nonlinear control (which has not been achieved so far) due to the 
network dynamics' being fundamentally nonlinear.  

Is it possible to use linear controllability as a kind of centrality measure
for complex networks? The answer is ``it depends.'' An essential requirement
for such a measure is the ability to distinguish and rank the nodes in the
network according to some criteria. Intuitively, one would 
hope that the nodes in the minimum controller set may be special and bear 
importance relative to other nodes. However, as demonstrated in our work, in 
a complex mutualistic network, the minimum controller set can be anything but 
unique. For a network of reasonable size, there is typically a vast number of
equivalent configurations or realizations of the set, a fact that was seldom
stated or studied in the existing literature of linear controllability of
complex networks. We note that, besides the linear structural~\cite{LSB:2011} 
and exact~\cite{YZDWL:2013} controllability theories, there are alternative 
frameworks such as the energy or linear Gramian based 
controllability~\cite{GPCTAKMVMG:2015}. However, the Gramian matrix depends 
on the chosen minimum controller set and the control signal input matrix. 
Our finding that, for some networks, almost all nodes can be in some 
realizations of the minimum controller set with approximately equal 
probability makes it difficult to use or exploit linear controllability
as a centrality measure for nodal ranking, such as network $A$ in
Fig.~\ref{fig:4Real_MN}(e). However, for other networks, some nodes 
are always or never in a driver set, which give a distribution of  
nodes in the minimum controller set. The distribution with respect to
the topology of the network may be informative and characteristic of some
empirical contexts~\cite{RR:2014,CARRSA:2017}.

The type of nonlinear control exploited in this paper for comparison with 
linear controllability is controlled management of the aftermath of a tipping 
point transition to enable species recovery. While this is a special type of 
control, its merit is rooted in the feasibility to quantify and rank the 
ability of individual nodes to promote recovery of the nonlinear dynamical 
network, so that the node-based, nonlinear control importance can be 
meaningfully compared with the corresponding linear control importance. Is 
there a more general approach to nonlinear network control which can be used 
for comparison with linear network control? We do not have an answer at 
the present, as the collective behaviors of nonlinear dynamical networks
are extremely diverse, so are the possible control strategies~\cite{WC:2002,
LWC:2004,SBGC:2007,YCL:2009,FMKS:2013,MFKS:2013,WSHWWGL:2016,ZYA:2017,
KSS:2017b}. However, regardless of the type of nonlinear control, heterogeneity
in the nodal importance ranking can be anticipated in general, due to the 
interplay between the nonlinear nodal dynamics and network structure. In 
contrast, as demonstrated in this paper, nodal importance ranking associated 
with linear controllability of complex networks exhibits a kind of 
heterogeneity opposite to that with nonlinear control, rendering linear 
controllability not useful for nonlinear dynamical networks in general. 

\section*{Methods}

\subsection*{General principle}

To obtain a statistical description of the roles played by the individual 
nodes and compare the nodal importance for nonlinear and linear control, 
we seek real world systems that meet the following two criteria: (a) the
underlying dynamical network is fundamentally nonlinear, for which a detailed 
mathematical description of the model is available, and (b) there exists
an issue of practical significance, with which nonlinear control is feasible.  
We find that mutualistic networks with a Holling type of 
dynamics~\cite{Holling:1959,Holling:1973} in ecology~\cite{BJMO:2003,GJT:2011,
NJB:2013,LNSB:2014,RSB:2014,DB:2014,GPJBT:2017,JHSLGHL:2018} and gene
regulatory networks with Michaelis-Menten type of dynamics in systems 
biology~\cite{Alon:2006,BBILA:2006,GBB:2016} satisfy these two criteria, with 
respect to the significant and broadly interesting issue of controlling 
tipping points. 

\subsection*{Nonlinear dynamical networks}

We have performed calculations and analyses for a large number of real-world 
pollinator-plant mutualistic networks available from the Web of Life database 
(http://www.web-of-life.es), which were reconstructed from empirical data 
collected from different geographic regions across different continents and 
climatic zones. The results reported in the main text are from the following
five representative mutualistic networks: (a) network $A$ ($N_A=38$ and 
$N_P=11$ with the number of mutualistic links $L=106$) from empirical data 
from Tenerife, Canary Islands~\cite{DHO:2003}, (b) network $B$ ($N_A=79$, 
$N_P=25$, and $L=299$) from Bristol, England~\cite{Memmott:1999}, (c) network 
$C$ ($N_A=36$, $N_P=61$, and $L=178$) from Morant Point, 
Jamaica~\cite{Percival:1974}, (d) network $D$ ($N_A=51$, $N_P=17$, and $L=129$)
from Tenerife, Canary Islands, and (e) network $E$ ($N_A=55$, $N_P=29$, 
and $L=145$) from Garajonay, Gomera, Spain. 

As a concrete example of gene regulatory networks, we study the transcription 
network of \textit{S. cerevisiae} of $4441$ nodes, for the representative
parameter setting~\cite{BBILA:2006} $B=1$, $f=1$, and $h=2$. In spite of 
the large number of genes involved in the network, the giant connected 
component in which each node can reach and is reachable from others along
a directed path has $60$ nodes only, and the size of the component in which 
each and every node has at least one incoming connection is $81$.  

\subsection*{Nonlinear control importance ranking}

For convenience, here we use the term ``nonlinear control importance'' to mean 
the statistical characterization of the nodal importance when carrying out
a physically meaningful type of control of the nonlinear dynamical network. 
Especially, we focus on controlling tipping points in complex pollinator-plant 
mutualistic networks and gene regulatory networks. 

For the mutualistic networks, a typical scenario for a tipping point to occur
is when the average mutualistic strength $\gamma_0$ is decreased towards zero.
The tipping point occurs at a critical value $\gamma_0^c$, at which the 
abundances of all species decrease to near zero values. There is global 
extinction for $\gamma_0 \le \gamma_0^c$. When $\gamma_0$ is increased from 
a value in the extinction region (e.g., in an attempt to restore the species 
abundances through improvement of the environment), recovery is not possible 
without control. A realistic control strategy was articulated, in which the 
abundance of a single pollinator species is maintained at a constant value, 
say $A_S$, through external means such as human management. 
We have observed numerically that, in the presence of control, a full recovery 
of all species abundances can be achieved - the phenomenon of ``control 
enabled recovery.'' 
For the same value of the controlled species level $A_S$, the critical 
$\gamma_0$ value of the recovery point depends on the particular species 
(node) subject to control. A smaller recovery point in $\gamma_0$ thus 
indicates that the control is more effective, which is species dependent. The 
species, or nodes in the network, can then be ranked with respect to the 
control. This provides a way to define the nodal importance associated with 
control of the underlying nonlinear network. In particular, let $\gamma_c^i$ 
be the system recovery point when the $i$th pollinator is subject to control. 
Choosing each and every pollinator species in turn as the controlled species, 
we obtain a set of values of the recovery point: $\{\gamma_c^i \}^{N_A}_{i=1}$.
Let $\gamma_c^{max}$ and $\gamma_c^{min}$ be the maximum and minimum values of 
the set. The importance of the pollinator species $i$ associated with control
of the tipping point can then be defined as 
\begin{equation} \label{eq:def_R_N}
R_{NL}^i=\frac{\gamma_c^{max}-\gamma_c^i}{\gamma_c^{max}-\gamma_c^{min}},
\end{equation}
where $0 \le R_{NL}^i \le 1$ and the control is more effective or, 
equivalently, the node subject to the control is more ``important'' if 
its corresponding value of $R_{NL}^i$ is larger. 

For the gene regulatory network, decreasing the value 
of the bifurcation parameter $C$ from one will result in a tipping point at
which the activities of all genes suddenly collapse to near zero values. The
behavior of sudden extinction at the tipping point can be harnessed by 
maintaining the activity level of a single active gene, e.g., the most active 
gene. In particular, when such control is present, the genes
``die'' in a benign way in that the death occurs one after another as the 
value of $C$ approaches zero, effectively eliminating the tipping point. 
We also find that, without control, it is not 
possible to recover the gene activities by increasing the value of $C$, but
a full recovery can be achieved with control. When a different gene is chosen
as the controlled target, for the same level of maintained activity, the
recovery point on the $C$ axis, denoted as $C_c$, is different, which provides 
the base to rank the ``importance'' of the genes with respect to control of
the nonlinear network. A gene with a relatively smaller value of $C_c$ is 
more important, as control targeted at it is more effective to restore the 
gene activities in the network.  

Similar to our approach to ranking the control importance for the 
pollinator-plant mutualistic networks, we define the following importance 
measure for gene $i$:
\begin{equation} \label{eq:def_R_N_g}
R_{NL}^i=\frac{C_c^{max}-C_c^i}{C_c^{max}-C_c^{min}},
\end{equation}
Where $C_c^i$ is the critical expression level to recover the whole 
system when the gene is subject to control, $C_c^{max}$ and $C_c^{min}$ are 
the maximum and minimum values of the recovery point among all the genes in
the network. 

\begin{figure}[ht]
\centering
\includegraphics[width=\linewidth]{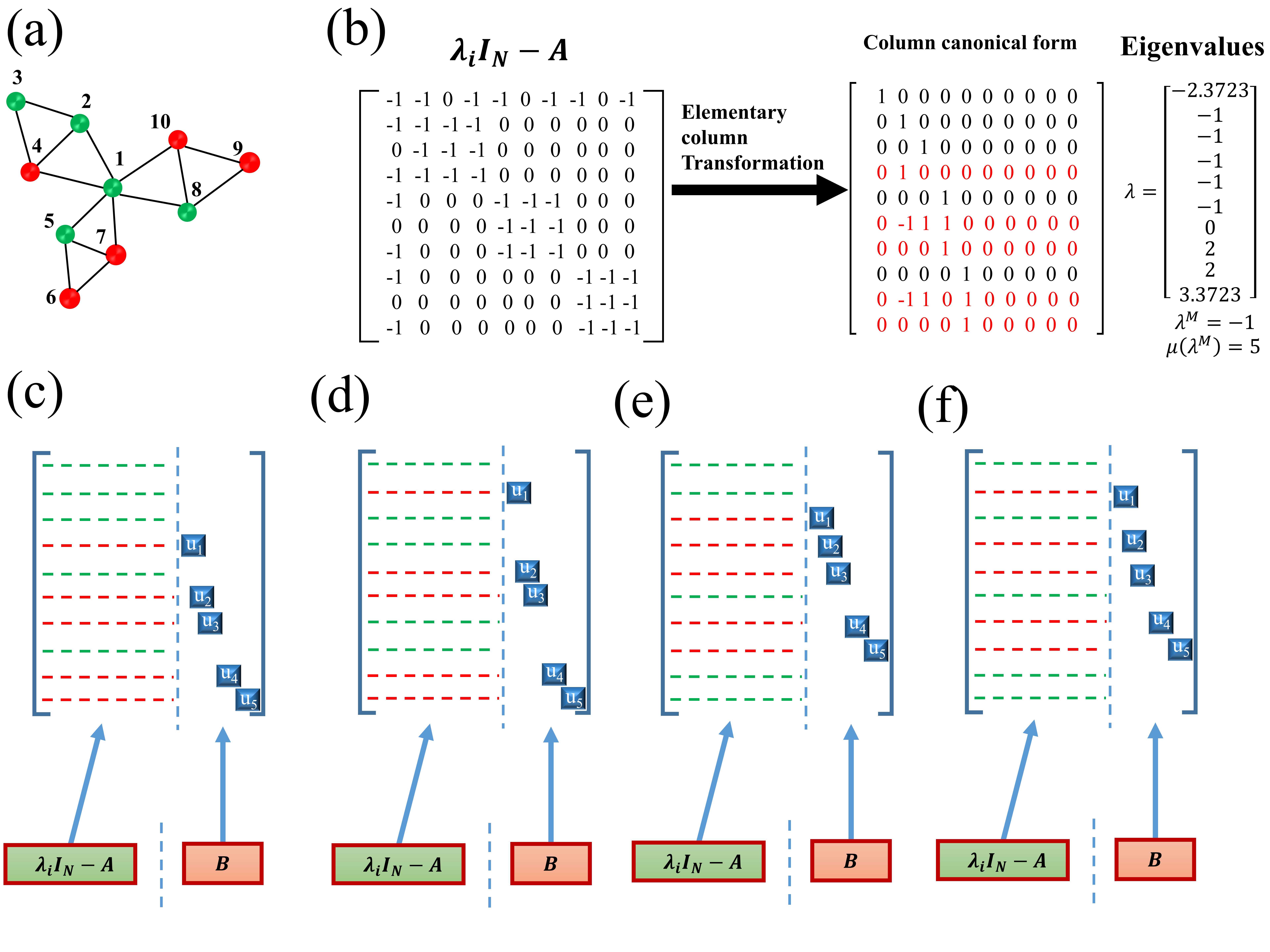}
\caption{ \textbf{Illustration of non-uniqueness of driver node set in linear
network control.} (a) A 10-node undirected network. Five eigenvalues of the
network connection matrix are identical: $\lambda = -1$, so its algebraic 
multiplicity is five. The geometric multiplicity of this eigenvalue is the  
number of linearly dependent rows in the matrix 
$\lambda_i \mathcal{I}_N - \mathcal{A}$,
which can be determined through elementary column transforms. (b) The matrix
$\lambda_i \mathcal{I}_N - \mathcal{A}$ and its representation after a 
series of elementary column transforms. The first, second, third, fifth, 
and eighth rows are distinct from all other rows, so they are linearly 
independent. The fourth, sixth, seventh, ninth and tenth rows are linearly 
dependent rows. Because there are five linearly independent and five linearly 
dependent rows, the number of ways to choose the latter is $54$. (c-f) Four 
distinct ways to choose the control input matrix to make the rank of the 
matrix $\lambda_i \mathcal{I}_N - \mathcal{A},\mathcal{B}$ ten. For this 
small network of size ten, any one of the nine out of the the ten nodes can 
be chosen to be a driver node.}
\label{fig:ECT_B}
\end{figure}

\subsection*{Linear control importance ranking}

Here, the term ``linear control importance ranking'' is referred to as the 
statistical ranking of the nodes in terms of their roles in the control
of the underlying linear dynamical network. This ranking can be determined by 
the exact linear controllability theory~\cite{YZDWL:2013}. To do so, we follow 
the existing studies that advocate the use of linear controllability 
for real world networked systems, such as those in 
Refs.~\cite{YVTCWSB:2017,GPCTAKMVMG:2015,MPGCGVB:2016,Tangetal:2017,TB:2018}. 
That is, we completely ignore the fact that the mutualistic network system 
and the gene regulatory network 
are highly nonlinear dynamical systems and instead treat them fictitiously as 
linear dynamical networks. For a network of $N$ nodes whose connecting topology
is characterized by the adjacent matrix $\mathcal{A}$, the linear control 
problem is formulated according to the following standard setting of 
canonically linear, time-invariant dynamical system:
\begin{equation} \label{eq:LSE}
\frac{d\mathbf{x}(t)}{dt}=\mathcal{A} \cdot \mathbf{x}(t) + 
\mathcal{B} \cdot \mathbf{u}(t),
\end{equation}
where $\mathbf{x}(t) \equiv (x_1(t),...,x_N(t))^T$ is the state vector of the 
system, $\mathcal{B}$ is the $N\times M$ input matrix ($M\leq N$) that 
specifies the control configuration - the set of $M$ nodes (driver nodes) to 
which external control signals $\mathbf{u}(t) = (u_1(t),...,u_M(t))^T$ 
should be applied. In general, the linear networked system Eq.~(\ref{eq:LSE}) 
can be controlled~\cite{Rugh:book} for properly chosen control vector 
$\mathbf{u}$ and for $M \ge N_D$, where $N_D$ is the minimum number of 
external signals required to fully control the network. The classic Kalman 
controllability rank condition~\cite{Kalman:1963} states that, system 
Eq.~(\ref{eq:LSE}) is controllable in the sense that it can be driven from 
any initial state to any desired final state in finite time if and only if 
the following $N\times NM$ controllability matrix
\begin{displaymath} 
\mathcal{C}=(\mathcal{B}, \mathcal{A}\cdot\mathcal{B}, \mathcal{A}^2\cdot\mathcal{B},\ldots,\mathcal{A}^{N-1}\cdot \mathcal{B}),
\end{displaymath}
has full rank:
\begin{displaymath} 
\mbox{rank}(C)=N.
\end{displaymath}
For a complex directed network, the linear structural controllability 
theory~\cite{Lin:1974} can be used to determine $N_D$ through identification
of maximum matching~\cite{LSB:2011}, the maximum set of links that do not 
share starting or ending nodes. A node is matched if there is a link in the 
maximum matching set points at it, and the directed network can be fully 
controlled if and only if there is a control signal on each unmatched node,
so $N_D$ is simply the number of unmatched nodes in the network. 

An alternative linear controllability framework, which is applicable to complex
networks of arbitrary topology (e.g., directed or undirected, weighted or 
unweighted), is the exact controllability theory~\cite{YZDWL:2013} derived 
from the PBH rank condition~\cite{PBH:1969}. In
particular, the linear system Eq.~(\ref{eq:LSE}) is fully controllable if and 
only if the following PBH rank condition 
\begin{equation} \label{eq:PBH_rank}
\mbox{rank}(c \mathcal{I}_N - \mathcal{A},\mathcal{B})=N,
\end{equation}
is met for any complex number $c$, where $\mathcal{I}_N$ is the $N\times N$ 
identity matrix. For any complex network defined by the general interaction 
matrix $\mathcal{A}$, it was proven~\cite{YZDWL:2013} that the network is 
fully controllable if and only if each and every eigenvalue $\lambda$ of 
$\mathcal{A}$ satisfies Eq.~(\ref{eq:PBH_rank}). For a set of control input 
matrices $\mathcal{B}$, $N_D$ can be determined as $N_D=min\{\mbox{rank}(B)\}$.
An equivalent but more practically useful criterion~\cite{YZDWL:2013} is that,
for a directed network, $N_D$ is nothing but the maximum geometric multiplicity
$\mu(\lambda_i)$ of the eigenvalue $\lambda_i$ of $A$:
\begin{equation} \label{eq:mini_set}
N_D = max_i\{\mu(\lambda_i)\},
\end{equation}
where $\lambda_i$ ($i=1,\ldots,l \le N$) are the distinct eigenvalues of 
$\mathcal{A}$ and geometric multiplicity of $\lambda_i$ is given by
\begin{displaymath}    
\mu(\lambda_i)=dimV_{\lambda_i} = N - \mbox{rank}(\lambda_i\mathcal{I}_N 
- \mathcal{A}).
\end{displaymath}
For a directed network, the exact controllability theory gives the same value
of $N_D$ as determined by the structural controllability theory. For an 
undirected network with arbitrary link weights, $N_D$ is determined by 
the maximum algebraic multiplicity (the eigenvalue degeneracy) 
$\delta(\lambda_i)$ of $\lambda_i$:
\begin{equation} \label{eq:ms_udnet}
N_D = max_i\{\delta(\lambda_i)\}.
\end{equation}
An issue of critical importance to our work but which is often ignored in the
existing literature on linear network controllability is the non-uniqueness of 
the set of the required driver nodes. In fact, for an arbitrary network with 
the value of $N_D$ determined, there can be a large number of equivalent 
configurations of the driver node set. This can be seen from the 
matrix $c \mathcal{I}_N - \mathcal{A}$ that appears in the PBH rank condition
Eq.~(\ref{eq:PBH_rank}). When $c$ is replaced by one of the eigenvalues of 
$\mathcal{A}$, say $\lambda_i$ (the one with the maximum algebraic 
multiplicity), the matrix $\lambda_i \mathcal{I}_N - \mathcal{A}$ contains
at least one dependent row. The quantity $N_D$ is nothing but the number of
linearly dependent rows of $\lambda_i \mathcal{I}_N - \mathcal{A}$. The 
control signals should then be applied to those nodes that correspond to 
the linearly dependent rows to make full rank the combined matrix 
$\lambda_i \mathcal{I}_N - \mathcal{A},\mathcal{B})$ in 
Eq.~(\ref{eq:PBH_rank}), as illustrated in Fig.~\ref{fig:ECT_B} for a small
network of size $N = 10$. The key fact is that there can be multiple but 
equivalent choices of the linearly dependent rows of the matrix 
$\lambda_i \mathcal{I}_N - \mathcal{A}$. For the small $10\times 10$ network 
in Fig.~\ref{fig:ECT_B}, there are $54$ such choices. The $N_D = 5$ driver 
nodes can then be chosen from the $N' = 9$ nodes as determined by the linearly 
dependent rows of $\lambda_i \mathcal{I}_N - \mathcal{A}$. When the network 
size $N$ is large, the $N_D \ll N$ driver nodes can be chosen from $N' \alt N$
nodes. Since $N_D \ll N'$, there can be great many distinct possibilities for
choosing the set of driver nodes (the number increases faster than exponential
with the network size). It is thus justified to define the probability for 
a node to be chosen as one of the driver nodes, so that the importance of 
each individual node in linear control can be determined. Specifically, the 
linear control importance of node $i$ can be defined as  
\begin{equation} \label{eq:def_R_L}
R_L^i = F_i/F,
\end{equation}    
where $F$ is the total number of configurations of the minimum controller sets
calculated and $F_i$ is the times that the $i^{th}$ node appears in these 
configurations. The probability $R_L^i$ thus gives the linear control 
importance ranking of the network, which can be meaningfully compared with 
the nonlinear control importance ranking.

\section*{Data Availability}

All relevant data are available from the authors upon request.

\section*{Code Availability}

All relevant computer codes are available from the authors upon request.

\section*{Acknowledgments}

We would like to acknowledge support from the Vannevar Bush Faculty
Fellowship program sponsored by the Basic Research Office of the Assistant
Secretary of Defense for Research and Engineering and funded by the Office
of Naval Research through Grant No.~N00014-16-1-2828.

\section*{Author Contributions}

YCL conceived the project. JJ performed computations and analysis. Both analyzed data. YCL wrote the paper with help from JJ.

\section*{Competing Interests}

The authors declare no competing interests.

\section*{Correspondence}

To whom correspondence should be addressed. E-mail: Ying-Cheng.Lai@asu.edu.

\appendix

\section{Linear controllability of \textit{C. elegans} connectome}

\begin{figure}[ht!]
\centering
\includegraphics[width=\linewidth]{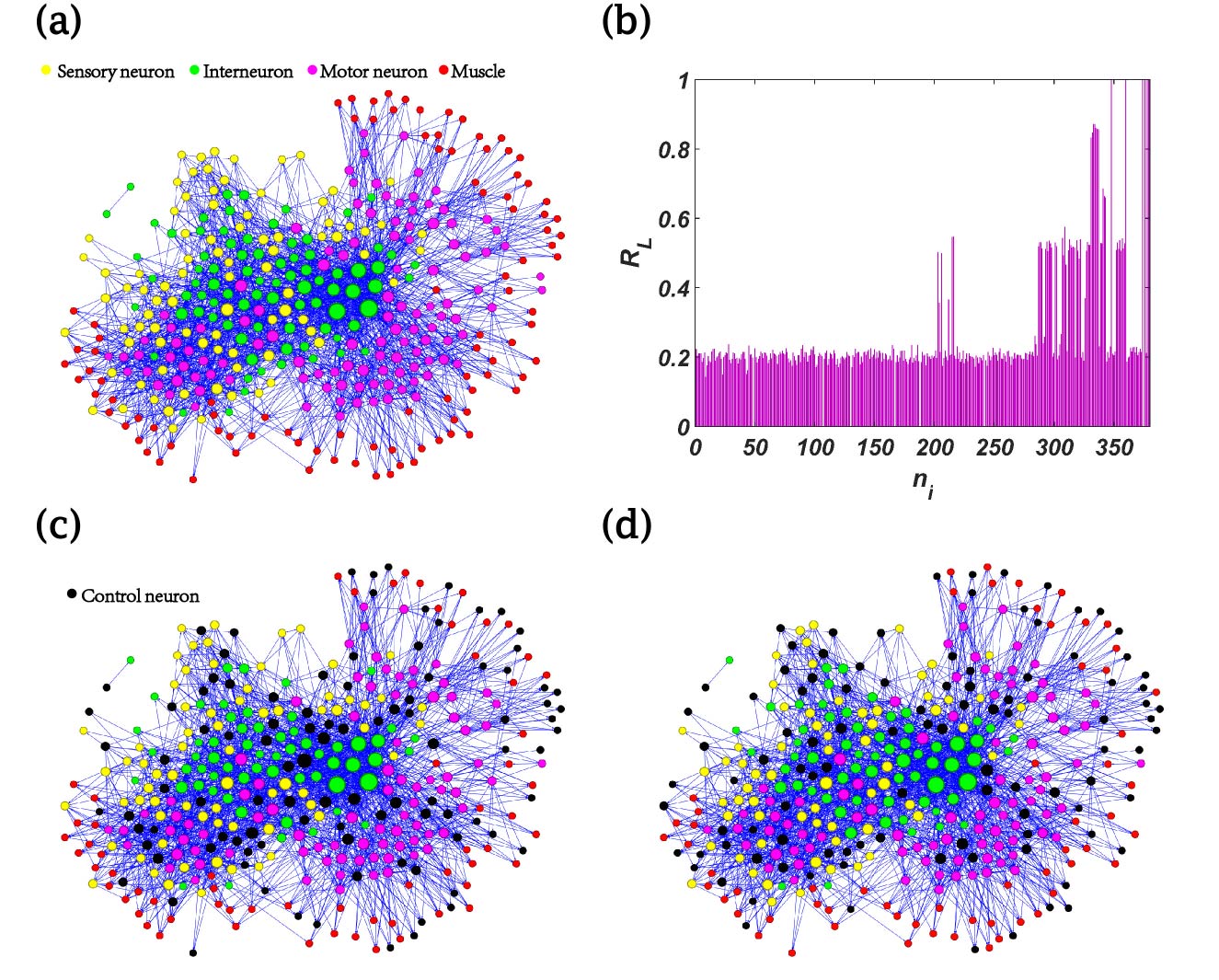}
\caption{ {\bf Linear control importance ranking for \textit{Caenorhabditis
elegans} connectome.} (a) A graphical representation of \textit{C. elegans}
connectome. The network contains $282$ neurons and $97$ muscle cells, where the
yellow, green, magenta, and red nodes represent sensory neurons, inter-neurons,
motor neurons, and muscle cells, respectively. The size of a node is
proportional to the sum of its in- and out-degrees. The dynamical network is
nonlinear, but a mathematical description of reasonable detail is not
available at the present. When the network is artificially treated as a
linear, time-invariant system, the size of the minimum controller set is
$N_D = 101$. (b) Linear control importance ranking, where the index on the
abscissa is arranged in a descending order of the nodal degree (sum of in-
and out-degrees). The importance distribution is approximately uniform across
the nodes, except for a small fraction of nodes. (c,d) Two realizations of
the minimum controller set (black dots). Because of the relatively large
size of the network, the number of distinct minimum controller sets is
quite large.}
\label{fig:LC_elegans}
\end{figure}

We report results from a linear controllability analysis of
\textit{C. elegans} connectome, whose network structure is shown in
Fig.~\ref{fig:LC_elegans}(a). In a recent work~\cite{YVTCWSB:2017}, the
neural network was treated as a linear, time invariant dynamical system
with control input signals applied to sensory neurons. It was found that
such a control signal would propagate to some motor neurons, and
the removal of one such neuron (that had not been identified previously)
would affect the muscle movement or function~\cite{YVTCWSB:2017}. We
have calculated that the size of the minimum controller set is quite large:
$N_D = 101$, which means that, since there are only $86$ sensory neurons
in \textit{C. elegans} connectome, it is not possible to control the linear
network even when each and every sensory neuron receives one independent
driving signal. There are many possible ways to place the required $N_D = 101$
control signals in the network, leading to many configurations of the
minimum controller set. We find that a typical realization of the set
contains both motor neurons and muscle cells. Figures~\ref{fig:LC_elegans}(b)
and \ref{fig:LC_elegans}(c) display two examples of the minimum controller
set, where the driver nodes are represented by black dots. The two
realizations share 43 common driver nodes, and the number of distinct drivers
is 58. Note the appearance of some muscle cells in both realizations.
Utilizing $1000$ random realizations, we calculate the linear control
importance ranking, as shown in Fig.~\ref{fig:LC_elegans}(d). It can be
seen that the statistical distribution of the importance is approximately
uniform for most nodes in the network, with only a few exceptions. There
is a probability for almost any neuron or muscle cell to belong to some
specific realization of the minimum controller set. We find that the average
values of the linear control importance for the three groups of neurons
(sensory, inter- and motor neurons) are approximately the same:
$\langle R_L\rangle_{SN}\approx 0.230$, $\langle R_L\rangle_{IN}\approx 0.211$,
and $\langle R_L \rangle_{MN} \approx 0.221$. However, the average linear nodal
importance for muscle cells is higher: $\langle R_L\rangle_{MC}\approx 0.399$.
These data indicate that the neurons in the connectome have equal chance to
be selected as a driver node, but a muscle cell is almost twice more likely
to appear in the minimum controller set. This result contradicts a general
understanding from both the biological and control perspectives, and has
intriguing implications to the relevance of the linear controllability
theory to \textit{C. elegans} connectome.
Specifically, from the point of view of biology, neurons send signals to the
muscle cells, but not the other way around. From the standpoint of actual
control of the network, a biologically meaningful driver set should favor
neurons. Yet the linear controllability theory gives the opposite result,
in contrast to the claim in Ref.~\cite{YVTCWSB:2017}.

\begin{figure}[ht]
\centering
\includegraphics[width = \linewidth]{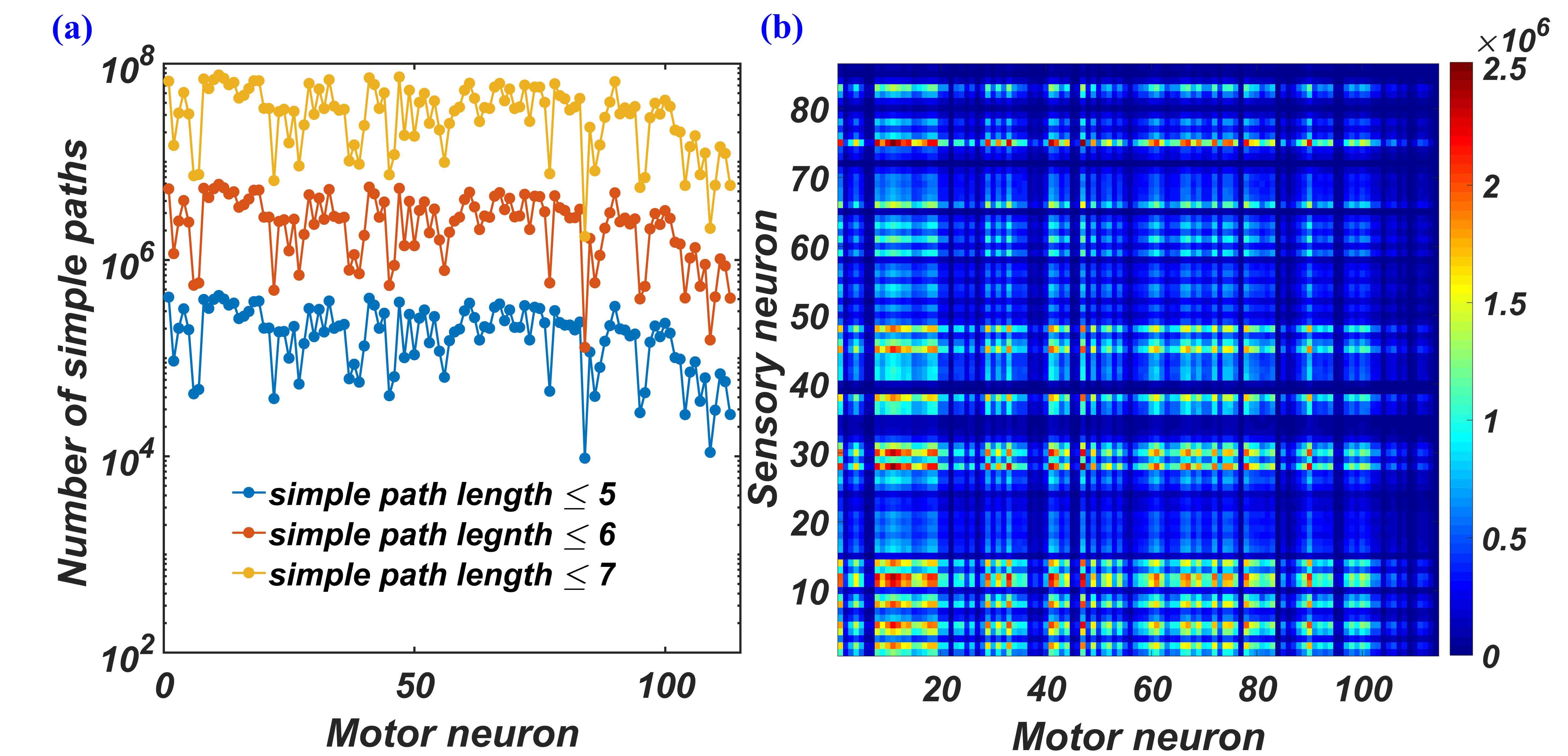}
\caption{ \textbf{Signal paths between sensory neurons and motor neurons in
\textit{C. elegans}}. (a) The numbers of direct paths from all sensory neurons
to motor neurons with path length less than or equal to seven (yellow),
six (red), and five (blue). The numbers of such paths are enormous. (b) Matrix
representation of all direct paths with path length less than or equal to
seven, from each and every sensory neuron to each and every motor neuron.}
\label{fig:C_elegans_SPC}
\end{figure}

In Ref.~\cite{YVTCWSB:2017}, some particular signal paths from the sensory
neurons to a special motor neuron were identified and deemed to be
particularly important based on the linear controllability theory.
Does linear control really reveal any specially important motor neurons, i.e.,
are there any differences among the motor neurons in terms of linear control
importance? To address this question, we map out all the direct paths among
the sensory and motor neurons that control the muscle cells and hence the
movement of \textit{C. elegans}. Figure~\ref{fig:C_elegans_SPC}(a) shows the
total numbers of direct paths of length $l$ less than or equal to five, six,
and seven from the sensory neurons to each and every motor neuron, where the
abscissa is the motor neuron index. The number of these paths is large. For
example, for $l \le 7$, for each and every motor neuron, there are between
$10^6$ and $10^8$ such paths. Apart from statistical fluctuations, the
numbers of paths are approximately constant across all the motor neurons,
suggesting the nonexistence of any special motor neuron. A matrix
representation of the paths for $l \le 7$ is shown in
Fig.~\ref{fig:C_elegans_SPC}(b). Between each and every pair of sensory
and motor neurons, the number of such paths is at least of the order of
$10^6$, although the numbers associated some specific paths can be about
two orders of magnitude higher.

\section{Description of real-world mutualistic networks}

\begin{table} [ht!]
        \centering
        \resizebox{1\textwidth}{!}{
                \begin{tabular}{|l|l|l|l|}
                        \hline
(1) Anastoechus, latifrons & (2) Anthophora, alluaudi & (3) Apis, mellifera & (4) Euodynerus, reflexus \\ \hline
(5) Geron, hesperidon & (6) Eristalis, tenax & (7) Megachile, canariensis & (8) Anthrax, anthrax \\ \hline
(9) Eucera, gracilipes & (10) Hyleaus, canariensis & (11) Lasioglossum, viride & (12) Linnaemyia, soror \\ \hline
(13) Cephalodromia & (14) Cyclyrius, webbianus & (15) Estheria, simonyi & (16) Lasioglossum, actifrons \\ \hline
(17) Melecta, curvispina & (18) Osmia, canariensis & 19) Andrena, wollestoni & (20) Colletes, dimidiatus \\ \hline
(21) Gasteruption  & (22) Lucilia, sericata & (23) Macroglossum, stellatarum & (24) Scaeva, albomaculata \\ \hline
(25) Stomorhina, lunata & (26) Unidentified & (27) Anthidium, manicatum & (28) Bibio, elmoi \\ \hline
(29) Dermasothes, gracile & (30) Drosophila & (31) Lasioglossum, chalcodes & (32) Leptochilus, eatoni \\ \hline
(33) Nyctia, lugubris & (34) Peleteria, ruficornis & (35) Phylloscopus, collybita & (36) Serinus, canarius \\ \hline
(37) Tachina, canariensis & (38) Tachysphex, unicolor &  & \\ \hline
                \end{tabular}
        }
        \caption{All names of the species
        in Fig.~1 in the main text.}
\end{table}

\begin{center}
\begin{longtable}{|c|c|c|c|p{0.3\linewidth}|}
\hline \multicolumn{1}{|c|}{\textbf{Index}} & \multicolumn{1}{c|}{\textbf{$\#$ Pollinators}} & \multicolumn{1}{c|}{\textbf{$\#$ Plants}} & \multicolumn{1}{c|}{\textbf{Linkage}} & \multicolumn{1}{c|}{\textbf{Network Location}} \\ \hline
\endfirsthead

\multicolumn{3}{c}%
{{\bfseries \tablename\ \thetable{} -- continued from previous page}} \\
                \hline \multicolumn{1}{|c|}{\textbf{Index}} &
                \multicolumn{1}{c|}{\textbf{$\#$ Pollinators}} &
                \multicolumn{1}{c|}{\textbf{$\#$ Plants}} &
                \multicolumn{1}{c|}{\textbf{Linkage}} &
                \multicolumn{1}{c|}{\textbf{Network Location}}\\ \hline
                \endhead

                1  &101 &84 &0.04 &Cordón del Cepo, Chile \tabularnewline \hline
                2  &64 &43 &0.07 &Cordón del Cepo, Chile \tabularnewline \hline
                3  &25 &36 &0.09 &Cordón del Cepo, Chile \tabularnewline \hline
                4  &102 &12 &0.14 &Central New Brunswick, Canada \tabularnewline \hline
                5  &275 &96 &0.03 &Pikes Peak, Colorado, USA \tabularnewline \hline
                6  &61  &17 &0.14 &Hickling, Norfolk, UK \tabularnewline \hline
                7  &36  &16 &0.15 &Shelfanger, Norfolk, UK \tabularnewline \hline
                8  &38  &11 &0.25 &Tenerife, Canary Islands \tabularnewline \hline
                9  &118 &24 &0.09 &Latnjajaure, Abisko, Sweden \tabularnewline \hline
                10 &76  &31 &0.19 &Zackenberg \tabularnewline \hline
                11 &13  &14 &0.29 &Mauritius Island \tabularnewline \hline
                12 &55  &29 &0.09 &Garajonay, Gomera, Spain \tabularnewline \hline
                13 &56  &9  &0.2  &KwaZulu-Natal region, South Africa \tabularnewline \hline
                14 &81  &29 &0.08 &Hazen Camp, Ellesmere Island, Canada \tabularnewline \hline
                15 &666 &131&0.03 &DaphnÃ­, Athens, Greece \tabularnewline \hline
16 &179 &26 &0.09 &Doñana National Park, Spain \tabularnewline \hline
                17 &79 &25 &0.15 &Bristol, England \tabularnewline \hline
                18 &108&36 &0.09 &Hestehaven, Denmark \tabularnewline \hline
                19 &85 &40 &0.08 &Snowy Mountains, Australia \tabularnewline \hline
                20 &91 &20 &0.1  &Hazen Camp, Ellesmere Island, Canada \tabularnewline \hline
                21 &677&91 &0.02 &Ashu, Kyoto, Japan \tabularnewline \hline
                22 &45 &21 &0.09 &Laguna Diamante, Mendoza, Argentina \tabularnewline \hline
                23 &72 &23 &0.08 &Rio Blanco, Mendoza, Argentina \tabularnewline \hline
                24 &18 &11 &0.19 &Melville Island, Canada \tabularnewline \hline
                25 &44 &13 &0.25 &North Carolina, USA \tabularnewline \hline
                26 &54 &105&0.04 &Galapagos \tabularnewline \hline
                27 &60 &18 &0.11 &Arthur's Pass, New Zealand \tabularnewline \hline
                28 &139&41 &0.07 &Cass, New Zealand \tabularnewline \hline
                29 &118&49 &0.06 &Craigieburn, New Zealand \tabularnewline \hline
                30 &53 &28 &0.07 &Guarico State, Venezuela \tabularnewline \hline
                31 &49 &48 &0.07 &Canaima Nat. Park, Venezuela \tabularnewline \hline
                32 &33 &7  &0.28 &Brownfield, Illinois, USA \tabularnewline \hline
33 &34 &13 &0.32 &Ottawa, Canada \tabularnewline \hline
                34 &128&26 &0.09 &Chiloe, Chile \tabularnewline \hline
                35 &36 &61 &0.08 &Morant Point, Jamaica \tabularnewline \hline
                36 &12 &10 &0.25 &Flores, AÃ§ores Island \tabularnewline \hline
                37 &40 &10 &0.18 &Hestehaven, Denmark \tabularnewline \hline
                38 &42 &8  &0.24 &Hestehaven, Denmark \tabularnewline \hline
                39 &51 &17 &0.15 &Tenerife, Canary Islands \tabularnewline \hline
                40 &43 &29 &0.09 &Windsor, The Cockpit Country, Jamaica \tabularnewline \hline
                41 &43 &31 &0.11 &Syndicate, Dominica \tabularnewline \hline
                42 &6  &12 &0.35 &Puerto Villamil, Isabela Island, Galapagos \tabularnewline \hline
                43 &82 &28 &0.11 &Hestehaven, Denmark \\\hline
                44 &609&110&0.02 &Amami-Ohsima Island, Japan \\\hline
                45 &26 &17 &0.14 &Uummannaq Island, Greenland \\\hline
                46 &44 &16 &0.39 &Denmark \\\hline
                47 &186&19 &0.12 &Isenbjerg \\\hline
                48 &236&30 &0.09 &Denmark \\\hline
                49 &225&37 &0.07 &Denmark \\\hline
                50 &35 &14 &0.18 &Tenerife, Canary Islands \\\hline
                51 &90 &14 &0.13 &Nahuel Huapi National Park, Argentina \\\hline
                52 &39 &15 &0.16 &Tundra, Greenladn \\\hline
                53 &294&99 &0.02 &Mt. Yufu, Japan \\\hline
                54 &318&113&0.02 &Kyoto City, Japan \\\hline
                55 &195&64 &0.03 &Nakaikemi marsh, Fukui Prefecture, Japan \\\hline
                56 &365&91 &0.03 &Mt. Kushigata, Yamanashi Pref., Japan \\\hline
                57 &883&114&0.02 &Kibune, Kyoto, Japan \\\hline
                58 &81 &32 &0.12 &Parc Natural del Cap de Creus \\\hline
                59 &13 &13 &0.42 &Parque Nacional do Catimbau \\\hline

\caption{The $59$ real pollinator-plant networks are from
web-of-life (http://www.web-of-life.es). For each network, the
linkage is normalized with respect to the corresponding fully
connected (all-to-all) network for which the linkage is $100\%$.}
\end{longtable}
\end{center}

%\bibliographystyle{naturemag}
%\bibliography{LC_NLC}

\end{document}